\documentclass[2p]{article}
\usepackage{amssymb,amsmath,amsthm}
\usepackage{indentfirst}
\usepackage{exscale}
\usepackage{relsize}
\usepackage[numbers,sort&compress]{natbib}

\usepackage{geometry}
\usepackage{color}
\newcommand{\R}{\mathbb{R}}
\theoremstyle{definition}

\allowdisplaybreaks
\theoremstyle{remark}

\numberwithin{equation}{section}
\UseRawInputEncoding

\begin{document}
\title{\Large\bf{ Nontrivial solutions for a $(p,q)$-Kirchhoff type system with concave-convex nonlinearities on locally finite graphs}}
\date{}
\author {Zhangyi Yu$^{1}$, \ Junping Xie$^{2}$\footnote{Corresponding author, E-mail address: hnxiejunping@163.com} , \ Xingyong Zhang$^{1,3}$\\
{\footnotesize $^1$Faculty of Science, Kunming University of Science and Technology, Kunming, Yunnan, 650500, P.R. China.}\\
{\footnotesize $^2$Faculty of Transportation Engineering, Kunming University of Science and Technology,}\\
 {\footnotesize Kunming, Yunnan, 650500, P.R. China.}\\
{\footnotesize $^{3}$Research Center for Mathematics and Interdisciplinary Sciences, Kunming University of Science and Technology,}\\
 {\footnotesize Kunming, Yunnan, 650500, P.R. China.}\\}

 \date{}
 \maketitle

 \begin{center}
 \begin{minipage}{15cm}
 \par
 \small  {\bf Abstract:} By using the well-known mountain pass theorem and Ekeland's variational principle, we prove that there exist at least two fully-non-trivial solutions for a $(p,q)$-Kirchhoff elliptic system with the Dirichlet boundary conditions and perturbation terms on a locally weighted and connected finite graph $G=(V,E)$.
 We also present a necessary condition of the existence of semi-trivial solutions for the system. Moreover, by using Ekeland's variational principle and Clark's Theorem, respectively, we prove that the system has at least one or multiple semi-trivial solutions when the perturbation terms satisfy different assumptions. Finally, we present a nonexistence result of solutions.


 \end{minipage}
 \end{center}
  \allowdisplaybreaks
 \vskip2mm
 {\section{Introduction }}
\setcounter{equation}{0}
In this paper, we are interested in the existence of nontrivial solutions to the nonhomogeneous $(p,q)$-Kirchhoff type system with the Dirichlet boundary condition on a locally weighted and connected finite graph $G=(V,E)$:
 \begin{eqnarray}
\label{eq1}
\begin{cases}
  -M_1(\|\nabla u\|_p^p)\Delta_pu=\lambda_1 h_1(x)|u|^{r-2}u+\frac{\alpha}{\alpha+\beta}h_2(x)|u|^{\alpha-2}u|v|^{\beta}+g_1(x),\;\;\;\;\hfill x\in\Omega,\\
  -M_2(\|\nabla v\|_q^q)\Delta_qv=\lambda_2 h_3(x)|v|^{r-2}v+\frac{\beta}{\alpha+\beta}h_2(x)|u|^{\alpha}|v|^{\beta-2}v+g_2(x),\;\;\;\;\hfill x\in\Omega,\\
  u(x)=v(x)=0,\;\;\;\;\hfill x\in\partial\Omega,
\end{cases}
\end{eqnarray}
where $V$ denotes the vertexes sets, $E$ denotes the edges sets,  $\Omega\cup\partial\Omega\subset V$ is a bounded domain, $M_i(s)=a_i+b_is^k,\;a_i,b_i>0,i=1,2,\;k\ge0$, $g_i:\Omega\to\R(i=1,2)$ are continuous functions which may change sign on $\Omega$, $p,q,r>1$, $\lambda_1,\lambda_2,\alpha,\beta>0$, $r<\min\{p,q\}\le(k+1)\max\{p,q\}<\alpha+\beta$ and $h_i:\Omega\to\R^+,i=1,2,3$. More detailed concepts on locally finite graphs can be seen in section 2 below.
\par
In recent years, the investigation on solutions for the equation on locally weighted finite graphs or weighted finite graphs, which essentially are  of discrete structure, gradually becomes a hot topic because the special definition of the Laplacian operator on graphs and the structure of graphs cause some new phenomenons. For example, in \cite{Yamabe 2016}, Grigor'yan-Lin-Yang proved a second order Laplacian equation on finite graphs and locally finite graphs admits a nontrivial solution with the help of the mountain pass theorem when the nonlinear term admits supequadratic growth. They also studied the $p$-Laplacian equation and generalized poly-Laplacian equation on finite graphs and locally finite graphs and established some similar results.
In \cite{Lin2017}, Lin-Wu investigated the existence and nonexistence of global solutions for a semilinear heat equation on a finite or locally finite connected weighted graph. They concluded that for a graph satisfying certain conditions, there exists a non-negative global solution when the initial value is small enough.
In \cite{Shao2024}, Shao-Yang-Zhao investigated a logarithmic type Schr\"{o}dinger equation on a finite connected subset of a locally finite graph with the nonlinear term $|u|^{p-2}u\log{u^2}$. For $p>2$, by virtue of the Brouwer degree theory and mountain pass theorem, they obtained the existence of ground state solutions.
In \cite{Yang2024}, Yang-Zhao investigated the nonlinear Schr\"{o}dinger equation with an $L^2$ mass constraint on both finite and locally finite graphs and obtained that the equation admits a normalized solution by utilizing variational methods.
In \cite{Yu 2023}, Yu-Zhang-Xie-Zhang were concerned with  the existence of nontrivial solutions for the generalized poly-Laplacian system with Dirichlet boundary condition on a locally finite graph $G=(V,E)$.
They establised sufficient conditions for the existence of a nontrivial solution for the system by employing the mountain pass theorem, with the nonlinear term satisfying asymptotically-$p$-linear conditions. In \cite{Pang 2024}, Pang-Xie-Zhang investigated several systems related to generalized poly-Laplacian and $(p,q)$-Laplacian on both weighed finite and locally finite graphs, with or without Dirichlet boundary values. They utilized an abstract critical points theorem without compactness conditions and showed that these systems have infinitely many nontrivial solutions with unbounded norms when the parameters locate some well-determined range. In \cite{Yang}, Yang-Zhang investigated the existence and multiplicity of nontrivial solutions for a $(p,q)$-Laplacian coupled system with perturbations and two parameters on locally finite graph. Using Ekeland's variational principle, they proved that the system has at least one nontrivial solution when the nonlinear term satisfies the sub-$(p,q)$ conditions. Furthermore, they established that the system has at least one solution of positive energy and one solution of negative energy when the nonlinear term satisfies the super-$(p,q)$ condition. The condition is weaker than the Ambrosetti-Rabinowitz condition, and their proof relied on the mountain pass theorem and Ekeland's variational principle.
Especially, in \cite{Yang 2023},  by using the mountain pass theorem and Ekeland's variational principle, Yang-Zhang investigated the existence of two nontrivial solutions for the following generalized poly-Laplacian system involving concave-convex nonlinearities and parameters with Dirichlet boundary value condition on a locally finite graph:
\begin{eqnarray}
\label{ed1}
\begin{cases}
  \pounds_{m_1,p}u=\lambda_1 h_1(x)|u|^{\gamma_1-2}u+\frac{\alpha}{\alpha+\beta}c(x)|u|^{\alpha-2}u|v|^{\beta},&\;\;\;\;\hfill x\in \Omega,\\
   \pounds_{m_2,q}v=\lambda_2 h_2(x)|v|^{\gamma_2-2}v+\frac{\beta}{\alpha+\beta}c(x)|u|^{\alpha}|v|^{\beta-2}v,&\;\;\;\;\hfill x\in \Omega,\\
   u=v=0,&\;\;\;\;\hfill x\in \partial\Omega.
\end{cases}
\end{eqnarray}
They obtained that system has at least one nontrivial solution of positive energy and one  nontrivial  solution of negative energy, respectively, where
$\pounds_{m_i,l}, i=1,2, l=p,q$ are the generalized poly-Laplacian operators and are the generalization of $l(=p,q)$-Laplacian operators.
For more related studies on graphs, we refer to \cite{Shao2023,Han2020,Han2021,Grigor2017,Ou2023,Pang 2023,Pinamonti2022,Zhang2018,Zhang2019,Zhang 2022}.
\par
Recently, in \cite{Pan2023}, an interesting work due to Pan-Ji is to investigate  the existence and convergence of the least energy sign-changing solutions to the following nonlinear Kirchhoff equation:
$$
-\left(a+b\int_V|\nabla u|^2d\mu\right)\Delta u+h(x)u=f(u),\;\;\;x\in V,
$$
where $a,b>0$, $G=(V,E)$ is a locally finite graph, $h:V\to\R$ be a function and $f:\R\to\R$. They proved the existence of a least energy sign-changing solution $u_b$ of the above equation if $c(x)$ and $f$ satisfy certain assumptions, and showed the energy of $u_b$ is strictly larger than twice that of the least energy solutions by using the constrained variational method.

\par
In this paper, motivated by \cite{Yang 2023,Pan2023} and following the line of thought in \cite{Yang 2023}, we would like to extend some results in \cite{Yang 2023} to the nonhomogeneous $(p,q)$-Kirchhoff system (\ref{eq1}). Moreover, we also present some results about the existence of semi-trivial solutions and the nonexistence of solutions. The main tools are the Ekeland's variational principle,  mountain pass theorem and Clark's Theorem which are recalled in section 2 below.
\par
Throughout this paper, we make the following assumptions.\\
$(H_1)$\; $M_i(s)=a_i+b_is^k,i=1,2,\;a_1>C^p_{1,p}(\Omega),\;a_2>C^q_{1,q}(\Omega),\;b_i>0,\;k\ge 0,\; s\ge 0$, where $C_{1,p}(\Omega)$ and $C_{1,q}(\Omega)$ are embedding constants given in Lemma 2.1 below.\\
$(H_2)\;$
\begin{eqnarray*}
& &g_i^{\star}:=\min_{x\in\Omega}g_i(x)\le g_i(x)\le\max_{x\in\Omega}g_i(x):=G_i,i=1,2,\\
& &0<h_i^{\star}:=\min_{x\in\Omega}h_i(x)\le h_i(x)\le\max_{x\in\Omega}h_i(x):=H_i,i=1,2,3.
\end{eqnarray*}
$(H_3)$ \;$p,q,r>1,\;\lambda_1,\lambda_2,\alpha,\beta>0,\;r<\min\{p,q\}\le\max\{p,q\}\le(k+1)\max\{p,q\}<\alpha+\beta$;\\
$(H_4)\;$
\begin{eqnarray*}
\begin{cases}
0<\lambda_1<a_1C^{-p}_{1,p}(\Omega)-1,\\
0<\lambda_2<a_2C^{-q}_{1,q}(\Omega)-1,\\
M_1\le\frac{\alpha+\beta}{\max\{p,q\}(k+1)}M_2,
\end{cases}
\end{eqnarray*}
and
\begin{eqnarray*}
& &\max\left\{\frac{\lambda_1(p-r)}{pr},\frac{\lambda_2(q-r)}{qr}\right\}\left(\|h_1\|^{\frac{p}{p-r}}_{L^{\frac{p}{p-r}}(\Omega)}+\|h_3\|^{\frac{q}{q-r}}_{L^{\frac{q}{q-r}}(\Omega)}\right)+\max\left\{\frac{p-1}{p},\frac{q-1}{q}\right\}\left(\|g_1\|^{\frac{p}{p-r}}_{L^{\frac{p}{p-r}}(\Omega)}+\|g_2\|^{\frac{q}{q-r}}_{L^{\frac{q}{q-r}}(\Omega)}\right)\\
&<&\frac{\alpha+\beta-\max\{p,q\}(k+1)}{\alpha+\beta}M^{\frac{\alpha+\beta}{\alpha+\beta-\max\{p,q\}(k+1)}}_1\left(\frac{\max\{p,q\}(k+1)}{(\alpha+\beta)M_2}\right)^{\frac{\max\{p,q\}(k+1)}{\alpha+\beta-\max\{p,q\}(k+1)}},
\end{eqnarray*}
where
\begin{eqnarray*}
& &M_1=2^{1-\max\{p,q\}(k+1)}\min\left\{\frac{b_1}{p(k+1)},\frac{b_2}{q(k+1)}\right\},\\
& &M_2=\frac{\|h_2\|_{\infty}}{(\alpha+\beta)^2}\max\left\{\alpha C^{\alpha+\beta}_{1,p}(\Omega),\beta C^{\alpha+\beta}_{1,q}(\Omega)\right\}.
\end{eqnarray*}

\par
Moreover, in this paper, we say that $(u,v)$ is a nontrivial solution of system $(\ref{eq1})$  if $(u,v)$ satisfies system $(\ref{eq1})$ and $(u,v)\not=(0,0)$, and we say that $(u,v)$ is a semi-trivial solution if $(u,v)$ satisfies system $(\ref{eq1})$, $(u,v)\not=(0,0)$ and one of $u$ and $v$ is zero. We say that $(u,v)$ is a fully-non-trivial solution if $(u,v)$ satisfies system $(\ref{eq1})$, $(u,v)\not=(0,0)$ and $(u,v)$ is not semi-trivial solution.
Our main results  read as follows.

\vskip2mm
\noindent
{\bf Theorem 1.1.} {\it Assume that $\Omega^{\circ}\neq\emptyset$, $\partial\Omega\neq\emptyset$ and $(H_1)$-$(H_3)$ hold true. If $\lambda_1$ and $\lambda_2$ satisfy $(H_4)$, system $(\ref{eq1})$ admits one nontrivial solution $(u_0,v_0)$ of positive energy and one nontrivial solution $(u^{\star}_0,v^{\star}_0)$ of negative energy.}

\vskip2mm
\par
The following Theorem 1.2 is a necessary condition of existence of semi-trivial solutions for system $(\ref{eq1})$.

\vskip2mm
\noindent
{\bf Theorem 1.2.} {\it  Assume that $\Omega^{\circ}\neq\emptyset$ and $\partial\Omega\neq\emptyset$. If $(u_0,0)$ is a semi-trivial solution of system $(\ref{eq1})$, then it holds that
$$
a_1\|u_0\|_{W^{1,p}_0(\Omega)}^p+b_1\|u_0\|_{W^{1,p}_0(\Omega)}^{p(k+1)}
\le\lambda_1 H_1 C_{1,p}^r(\Omega)\|u_0\|_{W^{1,p}_0(\Omega)}^{r}+G_1 C_{1,p}(\Omega)\|u_0\|_{W^{1,p}_0(\Omega)}.
$$
If $(0,v_0)$ is a semi-trivial solution of system $(\ref{eq1})$, then it holds that}
$$
a_2\|v_0\|_{W^{1,q}_0(\Omega)}^q+b_2\|v_0\|_{W^{1,q}_0(\Omega)}^{q(k+1)}
\le\lambda_2 H_2 C_{1,q}^r(\Omega)\|v_0\|_{W^{1,q}_0(\Omega)}^{r}+G_2 C_{1,q}(\Omega)\|v_0\|_{W^{1,q}_0(\Omega)}.
$$
\vskip2mm
\noindent
{\bf Remark 1.1.}  In \cite{Yang 2023}, Yang-Zhang did not present the results on the fully-non-trivial solutions. The main reason is that they assumed that $g_i\equiv0, i=1,2$ which make excluding the semi-trivial solutions become difficult. In our Theorem 1.1, it is easy to see that if $g_i\not\equiv 0,i=1,2$ hold, then both $(u_0,v_0)$ and $(u^{\star}_0,v^{\star}_0)$ are fully-non-trivial solutions. Hence, Theorem 1.1 and Theorem 1.2  can be seen as a generalization to Kirchoff type equation and completion of those results in \cite{Yang 2023} in some sense.

\vskip2mm
\par
The following Theorem 1.3 tells us that system $(\ref{eq1})$ has semi-trivial solutions.
\vskip2mm
\noindent
{\bf Theorem 1.3.} {\it Assume that $\Omega^{\circ}\neq\emptyset$ and $\partial\Omega\neq\emptyset$. If $g_1\not\equiv 0$ and  $g_2\equiv 0$ for all $x\in \Omega$, then system $(\ref{eq1})$ admits at least one semi-trivial solution $(u_0,0)$. If $g_2\not\equiv 0$ and  $g_1\equiv 0$ for all $x\in \Omega$, then system $(\ref{eq1})$ admits at least one semi-trivial solution $(0,v_0)$.  If $g_1\equiv 0$ and $g_2\equiv 0$ for all $x\in \Omega$, then for every $\lambda_1>0$, system $(\ref{eq1})$ has  $\#\Omega$  semi-trivial solutions  $\{(u_{0,j},0),j=1,2\cdots, \#\Omega\}$,  and for every $\lambda_2>0$, $\#\Omega$  semi-trivial solutions  $\{(0,v_{0,j}),j=1,2\cdots, \#\Omega\}$, where $\#\Omega$ denotes the number of elements in $\Omega$.}

\vskip2mm
\par
The following Theorem 1.4 is a nonexistence result of solutions for system $(\ref{eq1})$.
\vskip2mm
\noindent
{\bf Theorem 1.4.} {\it Suppose that $(\lambda_1,\lambda_2)\in \R^2$ satisfies that
\begin{eqnarray}
\label{eq38}
\lambda_1|s|^r\int_{\Omega}h_1(x)d\mu+\lambda_2|t|^r\int_{\Omega}h_3(x)d\mu+|s|^{\alpha}|t|^{\beta}\int_{\Omega}h_2(x)d\mu+s\int_{\Omega}g_1(x)d\mu+t\int_{\Omega}g_2(x)d\mu<0
\end{eqnarray}
for all $(s,t)\in \R^2$. Then system $(\ref{eq1})$ has no solutions.
}
\vskip2mm
{\section{Preliminaries}}
\setcounter{equation}{0}
\par
In this section, we review some basic concepts and properties on Sobolev spaces on locally finite graphs, which are taken from \cite{Yamabe 2016}.
\par
$G=(V, E)$ is called a locally connected finite graph if for any two vertexes in $V$ can be connected by finite edges in $E$ and for any $x\in V$ there are only finite edges $xy\in E$. The weight on any edge $xy\in E$ is represented by $\omega_{xy}$, which is required to satisfy $\omega_{xy}=\omega_{yx}$ and $\omega_{xy}>0$. We use $y\thicksim x$ to represent those vertices $y$ that connect with $x$. The distance  between two vertices $x$ and $y$, denoted by $d(x,y)$, is defined by the minimal number of edges which link $x$ to $y$. Assume $\Omega\subset V$. If $d(x,y)$ is uniformly bounded for any $x,y\in\Omega$, then $\Omega$ is known as a bounded domain in $V$. Let
$$
\partial\Omega=\{y\in V,\;y\notin\Omega|\;\exists\;x\in\Omega\;\text{such\;that}\;xy\in E\},
$$
which is known as the boundary of $\Omega$. Set $\Omega^{\circ}=\Omega\backslash\partial\Omega$ which is known as the interior of $\Omega$.
\par
Assume that $\mu:V\rightarrow \R^+$ is a finite measure  and there exists a $\mu_0>0$ such that $\mu(x)\ge \mu_0$.
The Laplacian operator $\Delta$ of $\psi$ is defined as
\begin{eqnarray}
\label{eq3}
\Delta \psi(x)=\frac{1}{\mu(x)}\sum\limits_{y\thicksim x}w_{xy}(\psi(y)-\psi(x)).
\end{eqnarray}
The corresponding gradient form has the expression given below
\begin{eqnarray}
\label{eq4}
\Gamma(\psi_1,\psi_2)(x)=\frac{1}{2\mu(x)}\sum\limits_{y\thicksim x}w_{xy}(\psi_1(y)-\psi_1(x))(\psi_2(y)-\psi_2(x)).
\end{eqnarray}
The length of the gradient is denoted by
\begin{eqnarray}
\label{eq5}
|\nabla \psi|(x)=\sqrt{\Gamma(\psi,\psi)(x)}=\left(\frac{1}{2\mu(x)}\sum\limits_{y\thicksim x}w_{xy}(\psi(y)-\psi(x))^2\right)^{\frac{1}{2}}.
\end{eqnarray}
For every fixed function $\psi:\Omega\rightarrow\mathbb{R}$, we set
\begin{eqnarray}
\label{eq7}
\int_\Omega \psi(x) d\mu=\sum\limits_{x\in \Omega}\mu(x)\psi(x).
\end{eqnarray}
\par
For any $l>1$, let
\begin{eqnarray}
\label{eq8}
\Delta_l\psi(x)=\frac{1}{2\mu(x)}\sum\limits_{y\sim x}\left(|\nabla \psi|^{l-2}(y)+|\nabla \psi|^{l-2}(x)\right)\omega_{xy}(\psi(y)-\psi(x)),
\end{eqnarray}
which is known as the $l$-Laplacian operator of $\psi$. Let
$$
C_c(\Omega)=\{\psi:V\rightarrow\R|supp\  \psi\subset\Omega\;\text{and}\;\psi(x)=0 \mbox{ for } \forall x\in V\backslash\Omega\}.
$$
For any function $\phi\in C_c(\Omega)$, there holds
\begin{eqnarray}
\label{eq15}
\int_{\Omega}\phi \Delta_l \psi d\mu=-\int_{\Omega\cup\partial\Omega}|\nabla \psi|^{l-2}\Gamma(\psi,\phi) d\mu.
\end{eqnarray}
\vskip2mm
\par
Assume that $G=(V,E)$ is a locally weighted finite graph and $\Omega\subset G$ is a bounded domain.
Let $1\le \gamma<+\infty$. Let $L^\gamma(\Omega)$ be the completion of $C_c(\Omega)$ under the norm
$$
\|u\|_{\gamma}=\left(\int_{\Omega}|u(x)|^\gamma d\mu\right)^{\frac{1}{\gamma}}.
$$
Let $W_0^{1,l}(\Omega)$ be the completion of $C_c(\Omega)$ under the norm
$$
\|u\|_{W_0^{1,l}(\Omega)}=\left(\int_{\Omega\cup\partial\Omega}|\nabla u(x)|^ld\mu\right)^{\frac{1}{l}}.
$$
 For any $u\in W_0^{1,l}(\Omega)$, we also define the following norm:
$$
\|u\|_{\infty,\Omega}=\max\limits_{x\in\Omega}|u(x)|.
$$

\vskip2mm
\noindent
{\bf Remark 2.1.} It is not difficult to obtain that  $W_0^{1,l}(\Omega)$ is of finite dimension and the dimension of $W_0^{1,l}(\Omega)$ is $\#\Omega$ which denotes the number of elements in $\Omega$. In fact, assume that $\Omega=\{x_1,\cdots,x_n\}$. If we let $u_i:V\to \R, i=1,\cdots,n$, and
$$
u_i(x)=
\begin{cases}
1,&x=x_i\\
0,&x\not=x_i,
\end{cases}\quad  i=1,\cdots,n.
$$
Then  $u_i\in C_c(\Omega), i=1,\cdots,n$ and  it is easy to verify that they are the basis of the space $W_0^{1,l}(\Omega)$.
 \vskip2mm
\noindent
{\bf Lemma 2.1.}  (\cite{Yamabe 2016}, Theorem 7 with $m=1$) {\it Suppose that $G=(V,E)$ is a locally finite graph, $\Omega$ is a bounded domain of $V$ such that $\Omega^{\circ}\neq\emptyset$. Let  $l>1$. Then $W_0^{1,l}(\Omega)$ is embedded in $L^{\gamma}(\Omega)$ for all $1\leq \gamma\leq+\infty$. Particularly, there exists a constant $C_{1,l}(\Omega)>0$ depending only on  $l$ and $\Omega$ such that
\begin{eqnarray}
\label{eq16}
\left(\int_{\Omega}|u|^\gamma d\mu\right)^{\frac{1}{\gamma}}\leq C_{1,l}(\Omega)\left(\int_{\Omega\cup\partial\Omega}|\nabla u|^ld\mu\right)^{\frac{1}{l}},
\end{eqnarray}
\begin{eqnarray}
\label{eq17}
\|u\|_{\infty,\Omega}\le M_l\|u\|_{W_0^{1,l}(\Omega)},
\end{eqnarray}
for all $1\le \gamma\le +\infty$ and all $u\in W_0^{1,l}(\Omega)$, where $C_{1,l}(\Omega)=\frac{C}{\mu_{\min,\Omega}}(1+|\sum_{x\in\Omega}\mu(x)|)$ with $C$ satisfying $\|u\|_{L^l(\Omega)}\le C\|u\|_{W_0^{1,l}(\Omega)}$, $M_l=\frac{C}{\mu_{\min}^\frac{1}{l}}$ and $\mu_{\min,\Omega}=\min_{x\in\Omega}\mu(x)$. Moreover, $W_0^{1,l}(\Omega)$ is pre-compact, that is, if $\{u_n\}$ is bounded in $W_0^{1,l}(\Omega)$, then up to a subsequence, there exists some $u\in W_0^{1,l}(\Omega)$ such that $u_n\rightarrow u$ in $W_0^{1,l}(\Omega)$.}

\vskip2mm
\par
Let $X$ be a Banach space and  $\varphi \in C^{1}(X,\R)$. A sequence is known as Palais-Smale sequence of a functional $\varphi$ if $\varphi(u_n)$ is bounded and
$\varphi'(u_n)\rightarrow 0$. The functional $\varphi$ satisfies the Palais-Smale (PS) condition if any Palais-Smale sequence of $\varphi$ has a convergent subsequence.

\vskip2mm
\noindent
{\bf Lemma 2.2.} (Mountain pass theorem  \cite{Rabinowitz 1986}) {\it Let $X$ be a real Banach space and $\psi \in C^{1}(X,\R)$, $\psi(0)=0$
and $\psi$ satisfies the (PS)-condition. Suppose that $\psi$ satisfies the following conditions:\\
(i) there exists a constant $ \varrho>0$ such that $ \psi|_{\partial B_{\varrho}(0)}> 0 $, where $B_\varrho=\{a\in X:\|a\|_X<\varrho\}$;\\
(ii) there exists $a\in X\backslash \bar B_{\varrho} (0)$ such that $\psi(a)\leq 0 $.\\
Then $\psi$ has a critical value $c$ with
$$
 c:=\inf_{\gamma\in\Gamma}\max_{t\in[0,1]}\psi(\gamma(t)),
$$
where
$$
 \Gamma:=\{\gamma\in C([0,1],X):\gamma(0)=0,\gamma(1)=a\}.
$$

}
\vskip2mm
\noindent
{\bf Lemma 2.3.} (A corollary of Ekeland's variational principle \cite{Mawhin2013}) {\it Assume that $\varphi\in C^1(X, \R)$ is bounded from below (above)
and satisfies the (PS) condition. Then
$c=\inf\limits_{u\in X} \varphi(u) (c=\sup\limits_{u\in X}
\varphi(u))$
is a critical value of $\varphi$.}

\vskip2mm
\noindent
{\bf Lemma 2.4.} (Clark's Theorem \cite{Rabinowitz 1986}) {\it  Let $X$ be a real Banach space and $\varphi$ be an even
function belonging to $C^1(X, \R)$ with $\varphi(0)=0$, bounded from below and satisfying
the (PS) condition. Denote by $\theta $ the zero element of $X$. $\Sigma$ indicates the family of sets $A \subset X \backslash \{\theta\}$
where $A$ is closed in $X$ and symmetric with respect to $\theta$, i.e. $u \in A$ implies $-u \in A$. Suppose that there is a set $K \in \Sigma$ such that $K$ is homeomorphic to $S^{j-1}$ ($j-1$ dimension unit sphere) by an odd map and $\sup\limits_K\varphi<0$. Then $\varphi$ has at
least $j$ distinct pairs of nonzero critical points.
 }

\vskip2mm
{\section{Proofs of main results}}
  \setcounter{equation}{0}

\vskip2mm
\par
In this section, we shall work in the space $X:=W_0^{1,p}(\Omega)\times W_0^{1,q}(\Omega)$ with the norm $\|(u,v)\|=\|u\|_{W_0^{1,p}(\Omega)}+\|v\|_{W_0^{1,q}(\Omega)}$. Then $X$ is a Banach space with finite dimension. Define the  corresponding variational functional $\varphi:X\to\R$ of system (\ref{eq1}) by
\begin{eqnarray}
\label{eq20}
\varphi(u,v)&=&\frac{a_1}{p}\|u\|_{W^{1,p}_0(\Omega)}^p+\frac{b_1}{p(k+1)}\|u\|_{W^{1,p}_0(\Omega)}^{p(k+1)}+\frac{a_2}{q}\|v\|_{W^{1,q}_0(\Omega)}^q+\frac{b_2}{q(k+1)}\|v\|_{W^{1,q}_0(\Omega)}^{q(k+1)}-\frac{\lambda_1}{r}\int_{\Omega}h_1(x)|u|^rd\mu\nonumber\\
& &-\frac{\lambda_2}{r}\int_{\Omega}h_3(x)|v|^rd\mu-\frac{1}{\alpha+\beta}\int_{\Omega}h_2(x)|u|^{\alpha}|v|^{\beta}d\mu-\int_{\Omega}g_1(x)ud\mu-\int_{\Omega}g_2(x)vd\mu.
\end{eqnarray}
Then a standard argument implies that  $\varphi(u,v)\in C^1(X,\R)$ (for example, see \cite{Yang}) and
\begin{eqnarray}
& &\langle\varphi'(u,v),(\phi_1,\phi_2)\rangle\nonumber\\
&=&\int_{\Omega\cup\partial\Omega}\left[(a_1+b_1\|u\|_{W^{1,p}_0(\Omega)}^{pk})|\nabla u|^{p-2}\nabla u\nabla\phi_1-\lambda_1 h_1(x)|u|^{r-2}u\phi_1-\frac{\alpha}{\alpha+\beta}h_2(x)|u|^{\alpha-2}u|v|^{\beta}\phi_1-g_1(x)\phi_1\right]d\mu\nonumber\\
& &+\int_{\Omega\cup\partial\Omega}\left[(a_2+b_2\|v\|_{W^{1,q}_0(\Omega)}^{qk})|\nabla v|^{q-2}\nabla v\nabla\phi_2-\lambda_2 h_3(x)|v|^{r-2}v\phi_2-\frac{\beta}{\alpha+\beta}h_2(x)|u|^{\alpha}|v|^{\beta-2}v\phi_2-g_2(x)\phi_2\right]d\mu\nonumber
\end{eqnarray}
for any $(u,v),(\phi_1,\phi_2)\in X$. The problem about finding the weak solutions for system (\ref{eq1}) is attributed to finding the critical points of the functional $\varphi$ on $X$. Actually, these weak solutions are also the classical solutions of system (\ref{eq1}) (the proof is similar to that in \cite{Zhang 2022}).

\vskip2mm
\noindent
{\bf Lemma 3.1.} {\it For every $(\lambda_1,\lambda_2)\in \R^2$ satisfying $(H_4)$, there exists a positive constant $\varrho_{(\lambda_1,\lambda_2)}$ such that $\varphi(u,v)>0$ if $\|(u,v)\|_X=\varrho_{(\lambda_1,\lambda_2)}$.
}
 \vskip0mm
 \noindent
{\bf Proof.}\; It follows from the Young's inequality and (\ref{eq16}) that
\begin{eqnarray}
\label{eq21}
\int_{\Omega}h_1(x)|u|^rd\mu&\le&\frac{p-r}{p}\int_{\Omega}h_1(x)^{\frac{p}{p-r}}d\mu+\frac{r}{p}\int_{\Omega}|u|^pd\mu\nonumber\\
&\le&\frac{p-r}{p}\|h_1\|_{L^{\frac{p}{p-r}}(\Omega)}^{\frac{p}{p-r}}+\frac{r}{p}C_{1,p}^p(\Omega)\|u\|_{W_0^{1,p}(\Omega)}^p.
\end{eqnarray}
Similarly,
\begin{eqnarray}
\label{eq22}
\int_{\Omega}h_3(x)|v|^rd\mu\le\frac{q-r}{q}\|h_3\|_{L^{\frac{q}{q-r}}(\Omega)}^{\frac{q}{q-r}}+\frac{r}{q}C_{1,q}^q(\Omega)\|v\|_{W_0^{1,q}(\Omega)}^q.
\end{eqnarray}
Moreover, we also have
\begin{eqnarray}
\label{eq23}
\int_{\Omega}|g_1||u|d\mu\le\frac{p-1}{p}\|g_1\|_{L^{\frac{p}{p-1}}(\Omega)}^{\frac{p}{p-1}}+\frac{1}{p}C_{1,p}^p(\Omega)\|u\|_{W_0^{1,p}(\Omega)}^p
\end{eqnarray}
and
\begin{eqnarray}
\label{eq24}
\int_{\Omega}|g_2||v|d\mu\le\frac{q-1}{q}\|g_2\|_{L^{\frac{q}{q-1}}(\Omega)}^{\frac{q}{q-1}}+\frac{1}{q}C_{1,q}^q(\Omega)\|v\|_{W_0^{1,q}(\Omega)}^q.
\end{eqnarray}
It follows from  Young's inequality and (\ref{eq16}) that
\begin{eqnarray}
\label{eq25}
\int_{\Omega}|h_2||u|^{\alpha}|v|^{\beta}d\mu&\le&\|h_2\|_{\infty,\Omega}\int_{\Omega}|u|^{\alpha}|v|^{\beta}d\mu\nonumber\\
&\le&\|h_2\|_{\infty}\left(\frac{\alpha}{\alpha+\beta}\int_{\Omega}|u|^{\alpha+\beta}d\mu+\frac{\beta}{\alpha+\beta}\int_{\Omega}|v|^{\alpha+\beta}d\mu\right)\nonumber\\
&\le&\|h_2\|_{\infty}\left(\frac{\alpha C_{1,p}^{\alpha+\beta}(\Omega)}{\alpha+\beta}\|u\|_{W_0^{1,p}(\Omega)}^{\alpha+\beta}+\frac{\beta C_{1,q}^{\alpha+\beta}(\Omega)}{\alpha+\beta}\|v\|_{W_0^{1,q}(\Omega)}^{\alpha+\beta}\right)\nonumber\\
&\le&\frac{\|h_2\|_{\infty}}{\alpha+\beta}\max\left\{\alpha C_{1,p}^{\alpha+\beta}(\Omega),\beta C_{1,q}^{\alpha+\beta}(\Omega)\right\}\|(u,v)\|_X^{\alpha+\beta}.
\end{eqnarray}
Thus, (\ref{eq20})$-$(\ref{eq25}) imply that when $(\lambda_1,\lambda_2)\in (0,a_1C_{1,p}^{-p}(\Omega)-1)\times(0,a_2C_{1,q}^{-q}(\Omega)-1)$, for any $(u,v)\in X$ with $\|(u,v)\|_X\le 1$,
\begin{eqnarray}
\label{eq28}
\varphi(u,v)&=&\frac{a_1}{p}\|u\|_{W^{1,p}_0(\Omega)}^p+\frac{b_1}{p(k+1)}\|u\|_{W^{1,p}_0(\Omega)}^{p(k+1)}+\frac{a_2}{q}\|v\|_{W^{1,q}_0(\Omega)}^q+\frac{b_2}{q(k+1)}\|v\|_{W^{1,q}_0(\Omega)}^{q(k+1)}-\frac{\lambda_1}{r}\int_{\Omega}h_1(x)|u|^rd\mu\nonumber\\
& &-\frac{\lambda_2}{r}\int_{\Omega}h_3(x)|v|^rd\mu-\frac{1}{\alpha+\beta}\int_{\Omega}h_2(x)|u|^{\alpha}|v|^{\beta}d\mu-\int_{\Omega}g_1(x)ud\mu-\int_{\Omega}g_2(x)vd\mu\nonumber\\
&\ge&\frac{1}{p}\left[a_1-(\lambda_1+1)C^p_{1,p}(\Omega)\right]\|u\|^p_{W^{1,p}_0(\Omega)}-\frac{\lambda_1(p-r)}{pr}\|h_1\|^{\frac{p}{p-r}}_{L^{\frac{p}{p-r}}(\Omega)}+\frac{1}{q}\left[a_2-(\lambda_2+1)C^q_{1,q}(\Omega)\right]\|v\|^q_{W^{1,q}_0(\Omega)}\nonumber\\
& &-\frac{\lambda_2(q-r)}{qr}\|h_3\|^{\frac{q}{q-r}}_{L^{\frac{q}{q-r}}(\Omega)}-\frac{\|h_2\|_{\infty}}{(\alpha+\beta)^2}\max\left\{\alpha C^{\alpha+\beta}_{1,p}(\Omega),\beta C^{\alpha+\beta}_{1,q}(\Omega)\right\}\|(u,v)\|^{\alpha+\beta}_X\nonumber\\
& &+\frac{b_1}{p(k+1)}\|u\|_{W^{1,p}_0(\Omega)}^{p(k+1)}+\frac{b_2}{q(k+1)}\|v\|_{W^{1,q}_0(\Omega)}^{q(k+1)}-\frac{p-1}{p}\|g_1\|^{\frac{p}{p-1}}_{L^{\frac{p}{p-1}}(\Omega)}-\frac{q-1}{q}\|g_2\|^{\frac{q}{q-1}}_{L^{\frac{q}{q-1}}(\Omega)}\nonumber\\
&\ge&\min\left\{\frac{a_1-(\lambda_1+1)C^p_{1,p}(\Omega)}{p},\frac{a_2-(\lambda_2+1)C^q_{1,q}(\Omega)}{q}\right\}\left(\|u\|^{\max\{p,q\}}_{W^{1,p}_0(\Omega)}+\|v\|^{\max\{p,q\}}_{W^{1,q}_0(\Omega)}\right)\nonumber\\
& &+\min\left\{\frac{b_1}{p(k+1)},\frac{b_2}{q(k+1)}\right\}\left(\|u\|^{\max\{p,q\}(k+1)}_{W^{1,p}_0(\Omega)}+\|v\|^{\max\{p,q\}(k+1)}_{W^{1,q}_0(\Omega)}\right)\nonumber\\
& &-\frac{\|h_2\|_{\infty}}{(\alpha+\beta)^2}\max\left\{\alpha C^{\alpha+\beta}_{1,p}(\Omega),\beta C^{\alpha+\beta}_{1,q}(\Omega)\right\}\|(u,v)\|^{\alpha+\beta}_X-\frac{\lambda_1(p-r)}{pr}\|h_1\|^{\frac{p}{p-r}}_{L^{\frac{p}{p-r}}(\Omega)}\nonumber\\
& &-\frac{\lambda_2(q-r)}{qr}\|h_3\|^{\frac{q}{q-r}}_{L^{\frac{q}{q-r}}(\Omega)}-\frac{p-1}{p}\|g_1\|^{\frac{p}{p-1}}_{L^{\frac{p}{p-1}}(\Omega)}-\frac{q-1}{q}\|g_2\|^{\frac{q}{q-1}}_{L^{\frac{q}{q-1}}(\Omega)}\nonumber\\
&\ge&2^{1-\max\{p,q\}}\min\left\{\frac{a_1-(\lambda_1+1)C^p_{1,p}(\Omega)}{p},\frac{a_2-(\lambda_2+1)C^q_{1,q}(\Omega)}{q}\right\}\|(u,v)\|^{\max\{p,q\}}_X\nonumber\\
& &+2^{1-\max\{p,q\}(k+1)}\min\left\{\frac{b_1}{p(k+1)},\frac{b_2}{q(k+1)}\right\}\|(u,v)\|^{\max\{p,q\}(k+1)}_X\nonumber\\
& &-\frac{\|h_2\|_{\infty}}{(\alpha+\beta)^2}\max\left\{\alpha C^{\alpha+\beta}_{1,p}(\Omega),\beta C^{\alpha+\beta}_{1,q}(\Omega)\right\}\|(u,v)\|^{\alpha+\beta}_X\nonumber\\
& &-\max\left\{\frac{\lambda_1(p-r)}{pr},\frac{\lambda_2(q-r)}{qr}\right\}\left(\|h_1\|^{\frac{p}{p-r}}_{L^{\frac{p}{p-r}}(\Omega)}+\|h_3\|^{\frac{q}{q-r}}_{L^{\frac{q}{q-r}}(\Omega)}\right)\nonumber\\
& &-\max\left\{\frac{p-1}{p},\frac{q-1}{q}\right\}\left(\|g_1\|^{\frac{p}{p-r}}_{L^{\frac{p}{p-r}}(\Omega)}+\|g_2\|^{\frac{q}{q-r}}_{L^{\frac{q}{q-r}}(\Omega)}\right)\nonumber\\
&\ge&2^{1-\max\{p,q\}(k+1)}\min\left\{\frac{b_1}{p(k+1)},\frac{b_2}{q(k+1)}\right\}\|(u,v)\|^{\max\{p,q\}(k+1)}_X\nonumber\\
& &-\frac{\|h_2\|_{\infty}}{(\alpha+\beta)^2}\max\left\{\alpha C^{\alpha+\beta}_{1,p}(\Omega),\beta C^{\alpha+\beta}_{1,q}(\Omega)\right\}\|(u,v)\|^{\alpha+\beta}_X\nonumber\\
& &-\max\left\{\frac{\lambda_1(p-r)}{pr},\frac{\lambda_2(q-r)}{qr}\right\}\left(\|h_1\|^{\frac{p}{p-r}}_{L^{\frac{p}{p-r}}(\Omega)}+\|h_3\|^{\frac{q}{q-r}}_{L^{\frac{q}{q-r}}(\Omega)}\right)\nonumber\\
& &-\max\left\{\frac{p-1}{p},\frac{q-1}{q}\right\}\left(\|g_1\|^{\frac{p}{p-r}}_{L^{\frac{p}{p-r}}(\Omega)}+\|g_2\|^{\frac{q}{q-r}}_{L^{\frac{q}{q-r}}(\Omega)}\right).
\end{eqnarray}
Note that
\begin{eqnarray*}
& &M_1=2^{1-\max\{p,q\}(k+1)}\min\left\{\frac{b_1}{p(k+1)},\frac{b_2}{q(k+1)}\right\},\\
& &M_2=\frac{\|h_2\|_{\infty}}{(\alpha+\beta)^2}\max\left\{\alpha C^{\alpha+\beta}_{1,p}(\Omega),\beta C^{\alpha+\beta}_{1,q}(\Omega)\right\}.
\end{eqnarray*}
Set
\begin{eqnarray}
\label{eq29}
g(z)&=&M_1z^{\max\{p,q\}(k+1)}-M_2z^{\alpha+\beta}\nonumber-\max\left\{\frac{\lambda_1(p-r)}{pr},\frac{\lambda_2(q-r)}{qr}\right\}\left(\|h_1\|^{\frac{p}{p-r}}_{L^{\frac{p}{p-r}}(\Omega)}+\|h_3\|^{\frac{q}{q-r}}_{L^{\frac{q}{q-r}}(\Omega)}\right)\nonumber\\
& &-\max\left\{\frac{p-1}{p},\frac{q-1}{q}\right\}\left(\|g_1\|^{\frac{p}{p-r}}_{L^{\frac{p}{p-r}}(\Omega)}+\|g_2\|^{\frac{q}{q-r}}_{L^{\frac{q}{q-r}}(\Omega)}\right),\;\;\;z\in[0,\infty).
\end{eqnarray}
To obtain $\varrho_{(\lambda_1,\lambda_2)}$ such  that $\varphi(u,v)>0$ when $\|(u,v)\|_X=\varrho_{(\lambda_1,\lambda_2)}$, we only need to prove that there is $z^{\star}_{(\lambda_1,\lambda_2)}\in (0,1]$ such that that $g(z^{\star}_{(\lambda_1,\lambda_2)})>0$. Indeed, it follows from (\ref{eq29}) that
\begin{eqnarray}
& &g'(z)=\max\{p,q\}(k+1)M_1z^{\max\{p,q\}(k+1)-1}-(\alpha+\beta)M_2z^{\alpha+\beta-1},\nonumber\\
& &g''(z)=\max\{p,q\}(k+1)(\max\{p,q\}(k+1)-1)M_1z^{\max\{p,q\}(k+1)-2}-(\alpha+\beta)(\alpha+\beta-1)M_2z^{\alpha+\beta-2}.\nonumber
\end{eqnarray}
Let $g'(z)=0$. Then the unique stationary point of $g$ can be obtained as follows:
\begin{eqnarray}
z^{\star}_{(\lambda_1,\lambda_2)}=\left(\frac{\max\{p,q\}(k+1)M_1}{(\alpha+\beta)M_2}\right)^{\frac{1}{\alpha+\beta-\max\{p,q\}(k+1)}}.\nonumber
\end{eqnarray}
By $(H_4)$, we know that $0<z^{\star}_{(\lambda_1,\lambda_2)}\le 1$. Furthermore, if follows from  $(k+1)\max\{p,q\}<\alpha+\beta$ that
\begin{eqnarray}
& &g''(z^{\star}_{(\lambda_1,\lambda_2)})\nonumber\\
&=&\max\{p,q\}(k+1)\left(\max\{p,q\}(k+1)-1\right)M_1\left(\frac{\max\{p,q\}(k+1)M_1}{(\alpha+\beta)M_2}\right)^{\frac{\max\{p,q\}(k+1)-2}{\alpha+\beta-\max\{p,q\}(k+1)}}\nonumber\\
& &-(\alpha+\beta)(\alpha+\beta-1)M_2\left(\frac{\max\{p,q\}(k+1)M_1}{(\alpha+\beta)M_2}\right)^{\frac{\alpha+\beta-2}{\alpha+\beta-\max\{p,q\}(k+1)}}\nonumber\\
&=&M_2\left(\frac{\max\{p,q\}(k+1)M_1}{(\alpha+\beta)M_2}\right)^{\frac{\alpha+\beta-2}{\alpha+\beta-\max\{p,q\}(k+1)}}[(\max\{p,q\}(k+1)-1)(\alpha+\beta)-(\alpha+\beta)(\alpha+\beta-1)]\nonumber\\
&=&(\max\{p,q\}(k+1)-\alpha-\beta)(\alpha+\beta)M_2\left(\frac{\max\{p,q\}(k+1)M_1}{(\alpha+\beta)M_2}\right)^{\frac{\alpha+\beta-2}{\alpha+\beta-\max\{p,q\}(k+1)}}\nonumber\\
&<&0.\nonumber
\end{eqnarray}
Note that $g(0)<0$ and $g(z)\to -\infty$ as $z\to+\infty$. Hence, $z^{\star}_{(\lambda_1,\lambda_2)}$ is the maximum point of $g$ and via $(H_4)$, there holds
\begin{eqnarray}
\max_{t\in[0,1]}g(z)&=&g(z^{\star}_{(\lambda_1,\lambda_2)})\nonumber\\
&=&M_1\left(\frac{\max\{p,q\}(k+1)M_1}{(\alpha+\beta)M_2}\right)^{\frac{\max\{p,q\}(k+1)}{\alpha+\beta-\max\{p,q\}(k+1)}}\nonumber\\
& &-M_2\left(\frac{\max\{p,q\}(k+1)M_1}{(\alpha+\beta)M_2}\right)^{\frac{\alpha+\beta}{\alpha+\beta-\max\{p,q\}(k+1)}}\nonumber\\
& &-\max\left\{\frac{\lambda_1(p-r)}{pr},\frac{\lambda_2(q-r)}{qr}\right\}\left(\|h_1\|^{\frac{p}{p-r}}_{L^{\frac{p}{p-r}}(\Omega)}+\|h_3\|^{\frac{q}{q-r}}_{L^{\frac{q}{q-r}}(\Omega)}\right)\nonumber\\
& &-\max\left\{\frac{p-1}{p},\frac{q-1}{q}\right\}\left(\|g_1\|^{\frac{p}{p-r}}_{L^{\frac{p}{p-r}}(\Omega)}+\|g_2\|^{\frac{q}{q-r}}_{L^{\frac{q}{q-r}}(\Omega)}\right)\nonumber\\
&=&\frac{\alpha+\beta-\max\{p,q\}(k+1)}{\alpha+\beta}M^{\frac{\alpha+\beta}{\alpha+\beta-\max\{p,q\}(k+1)}}_1\left(\frac{\max\{p,q\}(k+1)}{(\alpha+\beta)M_2}\right)^{\frac{\max\{p,q\}(k+1)}{\alpha+\beta-\max\{p,q\}(k+1)}}\nonumber\\
& &-\max\left\{\frac{\lambda_1(p-r)}{pr},\frac{\lambda_2(q-r)}{qr}\right\}\left(\|h_1\|^{\frac{p}{p-r}}_{L^{\frac{p}{p-r}}(\Omega)}+\|h_3\|^{\frac{q}{q-r}}_{L^{\frac{q}{q-r}}(\Omega)}\right)\nonumber\\
& &-\max\left\{\frac{p-1}{p},\frac{q-1}{q}\right\}\left(\|g_1\|^{\frac{p}{p-r}}_{L^{\frac{p}{p-r}}(\Omega)}+\|g_2\|^{\frac{q}{q-r}}_{L^{\frac{q}{q-r}}(\Omega)}\right)\nonumber\\
&>&0.\nonumber
\end{eqnarray}
If we let $\varrho_{(\lambda_1,\lambda_2)}=z^{\star}_{(\lambda_1,\lambda_2)}$, then the proof can be completed.
\vskip2mm
\noindent
{\bf Lemma 3.2.} {\it For every $(\lambda_1,\lambda_2)\in \R^2$ satisfying $(H_4)$, there exists $(u_{(\lambda_1,\lambda_2)},v_{(\lambda_1,\lambda_2)})\in X$ with $\|(u_{(\lambda_1,\lambda_2)},v_{(\lambda_1,\lambda_2)})\|_X>\varrho_{(\lambda_1,\lambda_2)}$ such that $\varphi(u_{(\lambda_1,\lambda_2)},v_{(\lambda_1,\lambda_2)})<0$.
}
 \vskip0mm
 \noindent
{\bf Proof.}\; For any given $(u,v)\in X$ satisfying $\int_{\Omega}h_2(x)|u|^{\alpha}|v|^{\beta}d\mu\not =0$ and any $t\in\R^+$, there holds
\begin{eqnarray}
\label{eq30}
& &\varphi(tu,tv)\nonumber\\
&=&\frac{a_1}{p}t^p\|u\|^p_{W^{1,p}_0(\Omega)}+\frac{b_1}{p(k+1)}t^{p(k+1)}\|u\|_{W^{1,p}_0(\Omega)}^{p(k+1)}+\frac{a_2}{q}t^q\|v\|_{W^{1,q}_0(\Omega)}^q+\frac{b_2}{q(k+1)}t^{q(k+1)}\|v\|_{W^{1,q}_0(\Omega)}^{q(k+1)}\nonumber\\
& &-\frac{\lambda_1}{r}t^r\int_{\Omega}h_1(x)|u|^rd\mu-\frac{\lambda_2}{r}t^r\int_{\Omega}h_3(x)|v|^rd\mu-\frac{1}{\alpha+\beta}t^{\alpha+\beta}\int_{\Omega}h_2(x)|u|^{\alpha}|v|^{\beta}d\mu-t\int_{\Omega}g_1(x)ud\mu-t\int_{\Omega}g_2(x)vd\mu\nonumber\\
&\le&\frac{a_1}{p}t^p\|u\|^p_{W^{1,p}_0(\Omega)}+\frac{b_1}{p(k+1)}t^{p(k+1)}\|u\|_{W^{1,p}_0(\Omega)}^{p(k+1)}+\frac{a_2}{q}t^q\|v\|_{W^{1,q}_0(\Omega)}^q+\frac{b_2}{q(k+1)}t^{q(k+1)}\|v\|_{W^{1,q}_0(\Omega)}^{q(k+1)}\nonumber\\
& &-\frac{\lambda_1h^{\star}_1}{r}t^r\int_{\Omega}|u|^rd\mu-\frac{\lambda_2h^{\star}_3}{r}t^r\int_{\Omega}|v|^rd\mu-\frac{h^{\star}_2}{\alpha+\beta}t^{\alpha+\beta}\int_{\Omega}|u|^{\alpha}|v|^{\beta}d\mu-tg^{\star}_1\int_{\Omega}ud\mu-tg^{\star}_2\int_{\Omega}vd\mu.
\end{eqnarray}
By virtue of the fact that $\alpha+\beta>(k+1)\max\{p,q\}$, it is easy  to obtain that there exists a sufficiently large $t_{(\lambda_1,\lambda_2)}$  satisfying $\|(t_{(\lambda_1,\lambda_2)}u,t_{(\lambda_1,\lambda_2)}v)\|_X>\varrho_{(\lambda_1,\lambda_2)}$ and $\varphi(t_{(\lambda_1,\lambda_2)}u,t_{(\lambda_1,\lambda_2)}v)<0$. If we choose $u_{(\lambda_1,\lambda_2)}=t_{(\lambda_1,\lambda_2)}u$ and $v_{(\lambda_1,\lambda_2)}=t_{(\lambda_1,\lambda_2)}v$, then the proof can be completed.
\qed
\vskip2mm
\noindent
{\bf Lemma 3.3.} {\it For every $(\lambda_1,\lambda_2)\in \R^2$ satisfying $(H_4)$,  the Palais-Smale condition holds for $\varphi$.
}
 \vskip0mm
 \noindent
{\bf Proof.}\; Assume that  $(u_m,v_m)\subseteq X$ is a Palais-Smale sequence, there exists a positive constant $C_0$  such that
$$
|\varphi(u_m,v_m)|\le C_0\;\;\text{for all}\;\;m\in\mathbb{N},\ \ \mbox{and }\varphi'(u_m,v_m)\to 0\;\;\text{as}\;\;m\to \infty.
$$
So,
\begin{eqnarray}
\label{eq31}
& &C_0+\|u_m\|_{W^{1,p}_0(\Omega)}+\|v_m\|_{W^{1,q}_0(\Omega)}\nonumber\\
&=&C_0+\|(u_m,v_m)\|_X\nonumber\\
&\ge&\varphi(u_m,v_m)-\frac{1}{\alpha+\beta}\langle\varphi'(u_m,v_m),(u_m,v_m)\rangle\nonumber\\
&=&\frac{a_1}{p}\|u_m\|_{W^{1,p}_0(\Omega)}^p+\frac{b_1}{p(k+1)}\|u_m\|_{W^{1,p}_0(\Omega)}^{p(k+1)}+\frac{a_2}{q}\|v_m\|_{W^{1,q}_0(\Omega)}^q+\frac{b_2}{q(k+1)}\|v_m\|_{W^{1,q}_0(\Omega)}^{q(k+1)}-\frac{\lambda_1}{r}\int_{\Omega}h_1(x)|u_m|^rd\mu\nonumber\\
& &-\frac{\lambda_2}{r}\int_{\Omega}h_3(x)|v_m|^rd\mu-\frac{1}{\alpha+\beta}\int_{\Omega}h_2(x)|u_m|^{\alpha}|v_m|^{\beta}d\mu-\int_{\Omega}g_1(x)u_md\mu-\int_{\Omega}g_2(x)v_md\mu\nonumber\\
& &-\frac{1}{\alpha+\beta}\bigg[\left(a_1+b_1\|u_m\|^{pk}_{W^{1,p}_0(\Omega)}\right)\|u_m\|^p_{W^{1,p}_0(\Omega)}-\lambda_1\int_{\Omega}h_1(x)|u_m|^rd\mu-\frac{\alpha}{\alpha+\beta}\int_{\Omega}h_2(x)|u_m|^{\alpha}|v_m|^{\beta}d\mu-\int_{\Omega}g_1(x)u_md\mu\nonumber\\
& &+\left(a_2+b_2\|v_m\|^{qk}_{W^{1,q}_0(\Omega)}\right)\|v_m\|^q_{W^{1,q}_0(\Omega)}-\lambda_2\int_{\Omega}h_3(x)|v_m|^rd\mu-\frac{\beta}{\alpha+\beta}\int_{\Omega}h_2(x)|u_m|^{\alpha}|v_m|^{\beta}d\mu-\int_{\Omega}g_2(x)v_md\mu\bigg]\nonumber\\
&=&a_1\left(\frac{1}{p}-\frac{1}{\alpha+\beta}\right)\|u_m\|^p_{W^{1,p}_0(\Omega)}+b_1\left(\frac{1}{p(k+1)}-\frac{1}{\alpha+\beta}\right)\|u_m\|^{p(k+1)}_{W^{1,p}_0(\Omega)}+a_2\left(\frac{1}{q}-\frac{1}{\alpha+\beta}\right)\|v_m\|^q_{W^{1,q}_0(\Omega)}\nonumber\\
& &+b_2\left(\frac{1}{q(k+1)}-\frac{1}{\alpha+\beta}\right)\|v_m\|^{q(k+1)}_{W^{1,q}_0(\Omega)}-\lambda_1\left(\frac{1}{r}-\frac{1}{\alpha+\beta}\right)\int_{\Omega}h_1(x)|u_m|^rd\mu-\lambda_2\left(\frac{1}{r}-\frac{1}{\alpha+\beta}\right)\int_{\Omega}h_3(x)|v_m|^rd\mu\nonumber\\
& &-\left(1-\frac{1}{\alpha+\beta}\right)\int_{\Omega}g_1(x)u_md\mu-\left(1-\frac{1}{\alpha+\beta}\right)\int_{\Omega}g_2(x)v_md\mu\nonumber\\
&\ge&a_1\left(\frac{1}{p}-\frac{1}{\alpha+\beta}\right)\|u_m\|^p_{W^{1,p}_0(\Omega)}+b_1\left(\frac{1}{p(k+1)}-\frac{1}{\alpha+\beta}\right)\|u_m\|^{p(k+1)}_{W^{1,p}_0(\Omega)}\nonumber\\
& &+a_2\left(\frac{1}{q}-\frac{1}{\alpha+\beta}\right)\|v_m\|^q_{W^{1,q}_0(\Omega)}+b_2\left(\frac{1}{q(k+1)}-\frac{1}{\alpha+\beta}\right)\|v_m\|^{q(k+1)}_{W^{1,q}_0(\Omega)}\nonumber\\
& &-\lambda_1\left(\frac{1}{r}-\frac{1}{\alpha+\beta}\right)H_1C^r_{1,p}(\Omega)\|u_m\|^r_{W^{1,p}_0(\Omega)}-\lambda_2\left(\frac{1}{r}-\frac{1}{\alpha+\beta}\right)H_3C^r_{1,q}(\Omega)\|v_m\|^r_{W^{1,q}_0(\Omega)}\nonumber\\
& &-\left(1-\frac{1}{\alpha+\beta}\right)\left(G_1C_{1,p}(\Omega)\|u_m\|_{W^{1,p}_0(\Omega)}+G_2C_{1,q}(\Omega)\|v_m\|_{W^{1,q}_0(\Omega)}\right).
\end{eqnarray}
Next, we prove that $\|(u_m,v_m)\|_X$ is bounded. Indeed, assume that
\begin{eqnarray}
\label{eq32}
\|u_m\|_{W^{1,p}_0(\Omega)}\to\infty\;\;\text{and}\;\;\|v_m\|_{W^{1,q}_0(\Omega)}\to\infty\;\;\text{as}\;\;m\to\infty.
\end{eqnarray}
Then for all large $m$, $(\ref{eq31})$ implies that
\begin{eqnarray*}
& &C_0+\|(u_m,v_m)\|_X\\
&\ge&\min\left\{a_1\left(\frac{1}{p}-\frac{1}{\alpha+\beta}\right),a_2\left(\frac{1}{q}-\frac{1}{\alpha+\beta}\right)\right\}\left(\|u_m\|^p_{W^{1,p}_0(\Omega)}+\|v_m\|^q_{W^{1,q}_0(\Omega)}\right)\\
& &+\min\left\{b_1\left(\frac{1}{p(k+1)}-\frac{1}{\alpha+\beta}\right),b_2\left(\frac{1}{q(k+1)}-\frac{1}{\alpha+\beta}\right)\right\}\left(\|u_m\|^{p(k+1)}_{W^{1,p}_0(\Omega)}+\|v_m\|^{q(k+1)}_{W^{1,q}_0(\Omega)}\right)\\
& &-\max\left\{\lambda_1\left(\frac{1}{r}-\frac{1}{\alpha+\beta}\right)H_1C^r_{1,p}(\Omega),\lambda_2\left(\frac{1}{r}-\frac{1}{\alpha+\beta}\right)H_3C^r_{1,q}(\Omega)\right\}\left(\|u_m\|^r_{W^{1,p}_0(\Omega)}+\|v_m\|^r_{W^{1,q}_0(\Omega)}\right)\\
& &-\max\left\{\left(1-\frac{1}{\alpha+\beta}\right)G_1C_{1,p}(\Omega),\left(1-\frac{1}{\alpha+\beta}\right)G_2C_{1,q}(\Omega)\right\}\left(\|u_m\|_{W^{1,p}_0(\Omega)}+\|v_m\|_{W^{1,q}_0(\Omega)}\right)\\
&\ge&2^{1-\min\{p,q\}}\min\left\{a_1\left(\frac{1}{p}-\frac{1}{\alpha+\beta}\right),a_2\left(\frac{1}{q}-\frac{1}{\alpha+\beta}\right)\right\}\|(u_m,v_m)\|^{\min\{p,q\}}_X\\
& &+2^{1-\min\{p,q\}(k+1)}\min\left\{b_1\left(\frac{1}{p(k+1)}-\frac{1}{\alpha+\beta}\right),b_2\left(\frac{1}{q(k+1)}-\frac{1}{\alpha+\beta}\right)\right\}\|(u_m,v_m)\|^{\min\{p,q\}(k+1)}_X\\
& &-\max\left\{\lambda_1\left(\frac{1}{r}-\frac{1}{\alpha+\beta}\right)H_1C^r_{1,p}(\Omega),\lambda_2\left(\frac{1}{r}-\frac{1}{\alpha+\beta}\right)H_3C^r_{1,q}(\Omega)\right\}\|(u_m,v_m)\|^r_X\\
& &-\max\left\{\left(1-\frac{1}{\alpha+\beta}\right)G_1C_{1,p}(\Omega),\left(1-\frac{1}{\alpha+\beta}\right)G_2C_{1,q}(\Omega)\right\}\|(u_m,v_m)\|_X,
\end{eqnarray*}
which implies that $\|(u_m,v_m)\|_X$ is bounded. Thus a contraction with $(\ref{eq32})$ appears.

\par
Assume that
\begin{eqnarray}
\label{eq33}
\|u_m\|_{W^{1,p}_0(\Omega)}\to\infty\;\;\text{as}\;\;m\to\infty,
\end{eqnarray}
and $\|v_m\|_{W^{1,q}_0(\Omega)}$ is bounded. Then it follows from $(\ref{eq31})$  that there exists a positive constant $C_1$  such that
\begin{eqnarray*}
& &C_1+\|u_m\|_{W^{1,p}_0(\Omega)}\\
&\ge&a_1\left(\frac{1}{p}-\frac{1}{\alpha+\beta}\right)\|u_m\|^p_{W^{1,p}_0(\Omega)}+b_1\left(\frac{1}{p(k+1)}-\frac{1}{\alpha+\beta}\right)\|u_m\|^{p(k+1)}_{W^{1,p}_0(\Omega)}\\
& &-\lambda_1\left(\frac{1}{r}-\frac{1}{\alpha+\beta}\right)H_1C^r_{1,p}(\Omega)\|u_m\|^r_{W^{1,p}_0(\Omega)}-\left(1-\frac{1}{\alpha+\beta}\right)G_1C_{1,p}(\Omega)\|u_m\|_{W^{1,p}_0(\Omega)},
\end{eqnarray*}
which implies that $\|u_m\|_X$ is bounded. Thus a contraction with $(\ref{eq33})$ appears.
 Similarly, assume that
$$
\|v_m\|_{W^{1,q}_0(\Omega)}\to\infty\;\;\text{as}\;\;m\to\infty,
$$
and $\|u_m\|_{W^{1,p}_0(\Omega)}$ is bounded,  a contradiction can also appear. Thus $\|u_m\|_{W^{1,p}_0(\Omega)}$ and $\|v_m\|_{W^{1,q}_0(\Omega)}$ are bounded. So the pre-compactness of $W_0^{1,l},l=p,q$ implies that there exist subsequences $\{u_{m_n}\}\subset\{u_m\}$ and $\{v_{m_n}\}\subset\{v_m\}$ such that $u_{m_n}\to u_0$ and $v_{m_n}\to v_0$ for some $u_0\in W^{1,p}_0(\Omega)$ and $v_0\in W^{1,q}_0(\Omega)$ as $n\to\infty$.
The proof is finished.
\qed

 \vskip2mm
 \noindent
{\bf Proof of Theorem 1.1.}\; It follows from  Lemmas 3.1-3.3 and Lemma 2.2 that for every $(\lambda_1,\lambda_2)\in \R^2$ satisfying $(H_4)$, system (\ref{eq1}) possesses  one nontrivial solution $(u_0,v_0)$ with positive energy. Moreover, note that $r<\min\{p,q\}$. Then, it follows from (\ref{eq30}) that there exists $t$ small enough such that
$$
\varphi(tu,tv)<0.
$$
So,
$$
-\infty<\inf\{\varphi(u,v):(u,v)\in\bar{B}_{\varrho_{(\lambda_1,\lambda_2)}}\}<0,
$$
where $\varrho_{(\lambda_1,\lambda_2)}$ is given in Lemma 3.1 and $\bar{B}_{\varrho_{(\lambda_1,\lambda_2)}}=\{(u,v)\in X:\|(u,v)\|_X\le \varrho_{(\lambda_1,\lambda_2)}\}$. Using the same proof as Theorem 1.3 in \cite{Yang 2023}, we can prove that system $(\ref{eq1})$ possesses a nontrivial solution $(u^{\star}_0,v^{\star}_0)$ with negative energy.
\qed

\vskip2mm
 \noindent
{\bf Proof of Theorem 1.2.}\; Suppose that $(u_0,0)$ is a semi-trivial solution of system $(\ref{eq1})$. Then
\begin{eqnarray}
0 &= &\langle\varphi'(u_0,0),(u_0,0)\rangle\nonumber\\
&=&\int_{\Omega\cup\partial\Omega}\left[(a_1+b_1\|u_0\|_{W^{1,p}_0(\Omega)}^{pk})|\nabla u_0|^{p}-\lambda_1 h_1(x)|u_0|^{r}-g_1(x)u_0\right]d\mu.\nonumber
\end{eqnarray}
Thus,
\begin{eqnarray}
    a_1\|u_0\|_{W^{1,p}_0(\Omega)}^{p}+b_1\|u_0\|_{W^{1,p}_0(\Omega)}^{p(k+1)}
&  =  &  \lambda_1 \int_{\Omega} h_1(x)|u_0|^{r}d\mu+\int_{\Omega}g_1(x)u_0d\mu\nonumber\\
& \le  & \lambda_1 H_1 C_{1,p}^r(\Omega)\|u_0\|_{W^{1,p}_0(\Omega)}^{r}+G_1 C_{1,p}(\Omega)\|u_0\|_{W^{1,p}_0(\Omega)}.
\end{eqnarray}
Similarly, suppose that $(0,v_0)$ is a semi-trivial solution of system $(\ref{eq1})$, it is easy to obtain that
$$
a_2\|v_0\|_{W^{1,q}_0(\Omega)}^q+b_2\|v_0\|_{W^{1,q}_0(\Omega)}^{q(k+1)}
\le\lambda_2 H_2 C_{1,q}(\Omega)^r\|v_0\|_{W^{1,q}_0(\Omega)}^{r}+G_2 C_{1,q}(\Omega)\|v_0\|_{W^{1,q}_0(\Omega)}.
$$
The proof is finished.
\qed

\vskip2mm
 \noindent
{\bf Proof of Theorem 1.3.}\; If  $g_1\not\equiv 0$ and  $g_2\equiv 0$. Then $v\equiv 0$ satisfies the second equation of  system $(\ref{eq1})$ and if we put it into the first equation of system  $(\ref{eq1})$. Then  system  $(\ref{eq1})$ reduces to the following equation:
 \begin{eqnarray}
\label{eeq1}
\begin{cases}
  -M_1(\|\nabla u\|_p^p)\Delta_pu=\lambda_1 h_1(x)|u|^{r-2}u+g_1(x),\;\;\;\;\hfill x\in\Omega,\\
  u(x)=0,\;\;\;\;\hfill x\in\partial\Omega.
\end{cases}
\end{eqnarray}
Next, we prove that equation (\ref{eeq1}) has at least one nontrivial solution. Indeed, a standard argument implies that the following corresponding variational functional $\psi$ of equation (\ref{eeq1}) is continuously differentiable on $W_0^{1,p}(\Omega)$, where
$$
\psi(u)=\frac{a_1}{p}\|u\|_{W^{1,p}_0(\Omega)}^p+\frac{b_1}{p(k+1)}\|u\|_{W^{1,p}_0(\Omega)}^{p(k+1)}-\frac{\lambda_1}{r}\int_{\Omega}h_1(x)|u|^rd\mu-\int_{\Omega}g_1(x)ud\mu.
$$
By Lemma 2.1, we have
\begin{eqnarray}\label{eeq3}
\psi(u)\ge \frac{a_1}{p}\|u\|_{W^{1,p}_0(\Omega)}^p+\frac{b_1}{p(k+1)}\|u\|_{W^{1,p}_0(\Omega)}^{p(k+1)}-\frac{\lambda_1}{r}H_1 C_{1,p}^r(\Omega)\|u\|_{W^{1,p}_0(\Omega)}^{r}-G_1 C_{1,p}(\Omega)\|u\|_{W^{1,p}_0(\Omega)}.
\end{eqnarray}
Since $r<p$, then $\psi$ is coercive and so it is bounded from blow. The coercivity of $\psi$ also implies that any Palais-Smale sequence $\{u_n\}$ is bounded. Furthermore, note that $W_0^{1,p}(\Omega)$ is pre-compact. Hence, $\{u_n\}$ has a convergent subsequence in $W_0^{1,p}(\Omega)$, still denoted by $\{u_n\}$, such that
$u_n\to u_0$ for some $u_0\in W_0^{1,p}(\Omega)$, which shows that $\psi$ satisfies the Palais-Smale condition. Thus, by Lemma 2.3, $u_0$ is a critical point of $\psi$ and so (\ref{eeq1}) has at least one solution $u_0$ in $W_0^{1,p}(\Omega)$. Obviously, $u_0\not=0$ since $g_1\not\equiv 0$.  It is easy to verify that $(u_0,0)$ satisfies system (\ref{eq1}).
 \par
 Similarly, if $g_2\not\equiv 0$ and  $g_1\equiv 0$.  We can also obtain that system (\ref{eq1}) has a semi-trivial solution $(0,v_0)$.
 \par
 If $g_1\equiv 0$ and  $g_2\equiv 0$ for all $x\in \Omega$,  then $v\equiv 0$ satisfies the second equation of  system $(\ref{eq1})$ and if we put it into the first equation of system  $(\ref{eq1})$. Then  system  $(\ref{eq1})$ reduces to the following equation
 \begin{eqnarray}
\label{eeq2}
\begin{cases}
  -M_1(\|\nabla u\|_p^p)\Delta_pu=\lambda_1 h_1(x)|u|^{r-2}u,\;\;\;\;\hfill x\in\Omega,\\
  u(x)=0,\;\;\;\;\hfill x\in\partial\Omega.
\end{cases}
\end{eqnarray}
A standard argument implies that the following corresponding variational functional $\psi$ of equation (\ref{eeq1}) is continuously differentiable on $W_0^{1,p}(\Omega)$, where
$$
\psi_0(u)=\frac{a_1}{p}\|u\|_{W^{1,p}_0(\Omega)}^p+\frac{b_1}{p(k+1)}\|u\|_{W^{1,p}_0(\Omega)}^{p(k+1)}-\frac{\lambda_1}{r}\int_{\Omega}h_1(x)|u|^rd\mu.
$$
Similar to the argument of (\ref{eeq3}), it is easy to obtain that $\psi_0$ is coercive, bounded from blow and satisfies the Palais-Smale condition. Moreover, it is obvious that $\psi_0$ is even on $W_0^{1,p}(\Omega)$ and $\psi_0(0)=0$. Since $W_0^{1,p}(\Omega)$ is of finite dimension, all the norms are equivalent. Hence, there exists a positive constant $D_0$ such that $\int_{\Omega}|u|^rd\mu\ge D_0\|u\|_{W^{1,p}_0(\Omega)}^r$ and then
$$
\psi_0(u)\le \frac{a_1}{p}\|u\|_{W^{1,p}_0(\Omega)}^p+\frac{b_1}{p(k+1)}\|u\|_{W^{1,p}_0(\Omega)}^{p(k+1)}-\frac{\lambda_1D_0h_1^{\star}}{r}\|u\|_{W^{1,p}_0(\Omega)}^r.
$$
For each $\lambda_1>0$, if we let
$$
K=\left\{u\in W_0^{1,p}(\Omega)\Big|\|u\|_{W^{1,p}_0(\Omega)}=d_{\lambda_1}\right\},
$$
then, on $K$ with choosing a sufficient small $d_{\lambda_1}$, we have
$$
\psi_0(u)\le \frac{a_1}{p}d^p+\frac{b_1}{p(k+1)}d^{p(k+1)}-\frac{\lambda_1D_0h_1^{\star}}{r}d^r<0,
$$
since $r<p$ and then $\sup\limits_K \psi_0<0$ since $K$ is closed. Moreover, it is easy to  obtain that $K$ is homeomorphic to $S^{\mbox{dim} W_0^{1,p}(\Omega)-1}$ ($\mbox{dim} W_0^{1,p}(\Omega)-1$ dimension unit sphere) by an odd map.
 Thus, $\psi_0$ satisfies all the conditions in Lemma 2.4. So $\psi_0$ has at least  $\mbox{dim} W_0^{1,p}(\Omega)$ nonzero critical points and then (\ref{eeq2}) has at least $\mbox{dim} W_0^{1,p}(\Omega)$ nonzero solutions $\{u_{0,j},j=1,2\cdots,\mbox{dim} W_0^{1,p}(\Omega)\}$, which shows that $\{(u_{0,j},0),j=1,2\cdots,\mbox{dim} W_0^{1,p}(\Omega)\}$  are semi-trivial solutions of system $(\ref{eq1})$.
 \par
 Similarly, similar to the above arguments, for each $\lambda_2>0$, it is not difficult to obtain that   $\{(0,v_{0,j}),j=1,2\cdots,\mbox{dim} W_0^{1,q}\}$  are also semi-trivial solutions of system $(\ref{eq1})$. Note that $\mbox{dim} W_0^{1,q}(\Omega)=\mbox{dim} W_0^{1,q}(\Omega)=\#\Omega$. Thus the proof is finished.
\qed

\vskip2mm
\noindent
{\bf Proof of Theorem 1.4.}\; Assume that $(\ref{eq1})$ has a solution $(\tilde{u},\tilde{v})$. Then $(\tilde{u},\tilde{v})$  is a critical point of $\varphi$  and
\begin{eqnarray*}
0&=&\langle\varphi'(\tilde{u},\tilde{v}),(\tilde{u},\tilde{v})\rangle\\
&=&\left(a_1+b_1\|\tilde{u}\|^{pk}_{W^{1,p}_0(\Omega)}\right)\|\tilde{u}\|^p_{W^{1,p}_0(\Omega)}-\lambda_1\int_{\Omega}h_1(x)|\tilde{u}|^rd\mu-\frac{\alpha}{\alpha+\beta}\int_{\Omega}h_2(x)|\tilde{u}|^{\alpha}|\tilde{v}|^{\beta}d\mu-\int_{\Omega}g_1(x)\tilde{u}d\mu\\
& &+\left(a_2+b_2\|\tilde{v}\|^{qk}_{W^{1,q}_0(\Omega)}\right)\|\tilde{v}\|^q_{W^{1,q}_0(\Omega)}-\lambda_2\int_{\Omega}h_3(x)|\tilde{v}|^rd\mu-\frac{\beta}{\alpha+\beta}\int_{\Omega}h_2(x)|\tilde{u}|^{\alpha}|\tilde{v}|^{\beta}d\mu-\int_{\Omega}g_2(x)\tilde{v}d\mu.
\end{eqnarray*}
So there holds
$$
\lambda_1\int_{\Omega}h_1(x)|\tilde{u}|^rd\mu+\lambda_2\int_{\Omega}h_3(x)|\tilde{v}|^rd\mu+\int_{\Omega}h_2(x)|\tilde{u}|^{\alpha}|\tilde{v}|^{\beta}d\mu+\int_{\Omega}g_1(x)\tilde{u}d\mu+\int_{\Omega}g_2(x)\tilde{v}d\mu\ge0.
$$
On the other hand, from $(\ref{eq38})$, we get
\begin{eqnarray*}
\lambda_1\int_{\Omega}h_1(x)|\tilde{u}|^rd\mu+\lambda_2\int_{\Omega}h_3(x)|\tilde{v}|^rd\mu+\int_{\Omega}h_2(x)|\tilde{u}|^{\alpha}|\tilde{v}|^{\beta}d\mu+\int_{\Omega}g_1(x)\tilde{u}d\mu+\int_{\Omega}g_2(x)\tilde{v}d\mu<0.
\end{eqnarray*}
Thus a contradiction appears. Hence,  system $(\ref{eq1})$ has no solutions.
\qed

\vskip2mm
\noindent
{\bf Remark 3.1.} Suppose that $g_i\equiv 0, i=1,2$. Then from Theorem 1.4 and its proofs, it is easy to obtain  that if
\begin{eqnarray*}
\lambda_1|s|^r\int_{\Omega}h_1(x)d\mu+\lambda_2|t|^r\int_{\Omega}h_3(x)d\mu+|s|^{\alpha}|t|^{\beta}\int_{\Omega}h_2(x)d\mu<0,
\end{eqnarray*}
then   system $(\ref{eq1})$ with $g_i\equiv 0, i=1,2$ has no nontrivial solutions.

 \vskip2mm
 \noindent
 {\bf Acknowledgments}\\
This work is supported by Yunnan Fundamental Research Projects of China (grant No: 202301AT070465) and  Xingdian Talent
Support Program for Young Talents of Yunnan Province in China.

\vskip2mm
\renewcommand\refname{References}
{}

\end{document}